\documentclass[10pt]{article}
\usepackage[cp866]{inputenc}
\begin{document}

\centerline {{\Large\bf Evolutionary forms: }}
\centerline {{\Large\bf The generation of differential-geometrical structures. }}
\centerline {{\Large\bf (Symmetries and Conservation laws.) }}

\centerline {L.I. Petrova}
\centerline{{\it Moscow State University, Russia, e-mail: ptr@cs.msu.su}}
\bigskip

Evolutionary forms, as well as exterior forms, are skew-symmetric
differential forms. But in contrast to the exterior forms, the basis of
evolutionary forms is deforming manifolds (manifolds with unclosed metric
forms). Such forms possess a peculiarity, namely, the closed inexact exterior 
forms are obtained from that.

The closure conditions of inexact exterior form (vanishing the differentials
of exterior and dual forms) point out to the fact that the closed inexact exterior
form is a quantity conserved on pseudostructure having the dual form as
the metric form. We obtain that the closed inexact exterior form and corresponding
dual form made up a conservative object, i.e. a quantity conserved on pseudostructure.
Such conservative object corresponds to the conservation law and is a
differential-geometrical structure.

Transition from the evolutionary form to the closed inexact exterior form
describes the process of generating the differential-geometrical structures.
This transition is possible only as a degenerate transformation, the condition of
which is a realization of a certain symmetry.

Physical structures that made up physical fields are such differential-geometrical
structures. And they are generated by material systems (medias). Relevant symmetries 
are caused by the degrees of freedom of material system.

\bigskip
{\bf Exterior and evolutionary forms.}

The difference between exterior and evolutionary skew-symmetric differential
forms is connected with the properties of manifolds
on which skew-symmetric forms are defined.

It is known that the exterior differential forms [1] are skew-symmetric
differential forms whose basis are differentiable manifolds or they can be
manifolds with structures of any type. Structures, on which exterior forms 
are defined, have closed metric forms.

It has been named as evolutionary forms skew-symmetrical differential forms
whose basis are deforming manifolds, i.e. manifolds with
unclosed metric forms. The metric form differential, and correspondingly 
its commutator, are nonzero. The commutators of metric forms of such 
manifolds describe a manifold deformation (torsion, curvature and so on).

Lagrangian manifolds,
manifolds constructed of trajectories of material system elements, tangent
manifolds of differential equations describing physical processes and others
can be examples of deforming manifolds.

A specific feature of the evolutionary forms, i.e skew-symmetric forms
defined on deforming manifolds, is the fact that commutators of these forms 
include commutators of the manifold metric forms being nonzero. Such commutators 
possess the evolutionary and topological properties. Just due to such properties 
of commutators the evolutionary forms can generate closed inexact exterior forms 
that correspond to the differential-geometrical structures.

It is known that the exterior differential form of degree $p$ ($p$-form)
can be written as [2,3]
$$
\theta^p=\sum_{i_1\dots i_p}a_{i_1\dots i_p}dx^{i_1}\wedge
dx^{i_2}\wedge\dots \wedge dx^{i_p}\quad 0\leq p\leq n\eqno(1)
$$
Here $a_{i_1\dots i_p}$ are the functions of the variables $x^{i_1}$,
$x^{i_2}$, \dots, $x^{i_p}$, $n$ is the dimension of space,
$\wedge$ is the operator of exterior multiplication, $dx^i$,
$dx^{i}\wedge dx^{j}$, $dx^{i}\wedge dx^{j}\wedge dx^{k}$, \dots\
is the local basis which satisfies the condition of exterior
multiplication.

The differential of the (exterior) form $\theta^p$ is expressed as
$$
d\theta^p=\sum_{i_1\dots i_p}da_{i_1\dots
i_p}\wedge dx^{i_1}\wedge dx^{i_2}\wedge \dots \wedge dx^{i_p} \eqno(2)
$$
The evolutionary differential form of degree $p$ ($p$-form) is written
similarly to exterior differential form. But the evolutionary form
differential cannot be written similarly to that presented for
exterior differential forms. In the evolutionary form
differential there appears an additional term connected with the fact
that the basis of the evolutionary form changes. For the differential 
forms defined on the manifold with unclosed metric form one has 
$d(dx^{\alpha_1}\wedge dx^{\alpha_2}\wedge \dots \wedge dx^{\alpha_p})\neq 0$.
(For the differential forms defined on the manifold with closed metric form 
one has $d(dx^{\alpha_1}\wedge dx^{\alpha_2}\wedge \dots \wedge dx^{\alpha_p})=0$). 
For this reason the differential of the evolutionary form $\theta^p$ can be
written as
$$
d\theta^p{=}\!\sum_{\alpha_1\dots\alpha_p}\!da_{\alpha_1\dots\alpha_p}\wedge
dx^{\alpha_1}\wedge dx^{\alpha_2}\dots \wedge
dx^{\alpha_p}{+}\!\sum_{\alpha_1\dots\alpha_p}\!a_{\alpha_1\dots\alpha_p}
d(dx^{\alpha_1}\wedge dx^{\alpha_2}\dots \wedge dx^{\alpha_p})\eqno(3)
$$
where the second term is a differential of unclosed metric form being
nonzero.

[In further presentation the symbol of summing $\sum$ and the symbol
of exterior multiplication $\wedge$ will be omitted. Summation
over repeated indices is implied.]

The second term connected with the differential of the basis can be expressed
in terms of the metric form commutator.

For example, let us consider the first-degree form
$\theta=a_\alpha dx^\alpha$. The differential of this form can
be written as
$$d\theta=K_{\alpha\beta}dx^\alpha dx^\beta\eqno(4)$$
where
$K_{\alpha\beta}=a_{\beta;\alpha}-a_{\alpha;\beta}$ are 
components of the commutator of the form $\theta$, and
$a_{\beta;\alpha}$, $a_{\alpha;\beta}$ are covariant
derivatives. If we express the covariant derivatives in terms of
the connectedness (if it is possible), they can be written
as $a_{\beta;\alpha}=\partial a_\beta/\partial
x^\alpha+\Gamma^\sigma_{\beta\alpha}a_\sigma$, where the first
term results from differentiating the form coefficients, and the
second term results from differentiating the basis. We arrive at the
following expression for the commutator components of the form $\theta$
$$
K_{\alpha\beta}=\left(\frac{\partial a_\beta}{\partial
x^\alpha}-\frac{\partial a_\alpha}{\partial
x^\beta}\right)+(\Gamma^\sigma_{\beta\alpha}-
\Gamma^\sigma_{\alpha\beta})a_\sigma\eqno(5)
$$
Here the expressions
$(\Gamma^\sigma_{\beta\alpha}-\Gamma^\sigma_{\alpha\beta})$
entered into the second term are just the components of
commutator of the first-degree metric form.

If to substitute the expressions (5) for evolutionary form
commutator into formula (4), we obtain the
following expression for the differential of the first degree
skew-symmetric form
$$
d\theta=\left(\frac{\partial a_\beta}{\partial
x^\alpha}-\frac{\partial a_\alpha}{\partial
x^\beta}\right)dx^\alpha dx^\beta+\left((\Gamma^\sigma_{\beta\alpha}-
\Gamma^\sigma_{\alpha\beta})a_\sigma\right)dx^\alpha dx^\beta\eqno(6)
$$
The second term in the expression for the differential of skew-symmetric form
is connected with the differential of the manifold metric form, which is expressed
in terms of the metric form commutator.

Thus, the differentials and, correspondingly, the commutators of
exterior and evolutionary forms are of different types. In
contrast to the exterior form commutator, the evolutionary form
commutator includes two terms. These two terms have different
nature, namely, one term is connected with coefficients of the
evolutionary form itself, and the other term is connected with
differential characteristics of manifold. Interaction between
terms of the evolutionary form commutator (interactions between
coefficients of evolutionary form and its basis) provide the
foundation of evolutionary processes that lead to generation of
closed inexact exterior forms to which the
differential-geometrical structures are assigned.

\bigskip
{\bf Closed inexact exterior forms. Conservation laws.

Differential-geometrical structures.}

From the closure condition of exterior form $\theta^p$:
$$
d\theta^p=0\eqno(7)
$$
one can see that the closed exterior form $\theta^p$ is a conserved
quantity. This means that it can correspond to a conservation law,
namely, to some conservative quantity.

If the form is closed only on pseudostructure, i.e. this form is a closed
inexact one, the closure conditions are written as
$$
d_\pi\theta^p=0\eqno(8)
$$
$$
d_\pi{}^*\theta^p=0\eqno(9)
$$
where ${}^*\theta^p$ is the dual form.

Condition (9), i.e. the closure condition for dual form, specifies
a pseudostructure $\pi$.
\{Cohomology (de Rham cohomology,
singular cohomology), sections of cotangent bundles and
so on may be regarded as examples of pseudostructures.\}
From conditions (8) and (9) one can see the following. The dual form
(pseudostructure) and closed inexact form (conservative quantity)
made up a conservative object that can also correspond to some
conservation law.
Conservative object, which corresponds to the conservation law,
is a differential-geometrical structure. Such
differential-geometrical structures are examples of G-structures.
The physical structures, which forms physical
fields, and corresponding conservation laws are just such structures.

The properties of closed exterior forms reflect also the properties of
differential-geometrical structures. 
The following properties should be emphasized.

1. {\it Invariance.}

It is known that the closed exact form is
a differential of the form of lower degree:
$$
\theta^p=d\theta^{p-1}\eqno(10)
$$
Closed inexact form is also a differential, and yet not a total one but
an interior on pseudostructure
$$
\theta^p_\pi=d_\pi\theta^{p-1}\eqno(11)
$$

Since the closed exterior differential forms are differentials,
they turn out to be invariant under all transformations that conserve the
differential. Gauge transformations (the unitary transformations,
tangent, canonical, and gradient ones) are examples of such
nondegenerate transformations under which closed exterior forms,
and hence, the differential-geometrical structures as well, turn out to
be invariant.

2. {\it Conjugacy}

Closure of exterior differential forms is a result of the conjugacy of
elements of exterior or dual forms. The closure property of the exterior
form means that any objects, namely, elements of the exterior form,
components of elements, elements of the form differential, exterior and
dual forms and others, turn out to be conjugated.

Conjugacy is possible if there is one or another type of symmetry.

Gauge symmetries are the symmetries of closed exterior differential
forms. They are obtained as the result of conjugacy of any
exterior form elements. The physical structures, which are
differential-geometrical structures and correspond to conservation laws,
are connected with gauge symmetries.

Mathematically these properties of closed exterior forms are written as
identical relations.
Since the conjugacy is a certain connection between two operators or
mathematical objects, it is evident that, to express a conjugacy
mathematically, it can be used relations.
These are identical relations.

The identical relations express the fact that each closed exterior
form is a differential of some exterior form (with a degree less
by one). In general form such an identical relation can be written as
$$
d _{\pi}\varphi=\theta _{\pi}^p\eqno(12)
$$

In this relation the form in the right-hand side has to be a
{\it closed} one.

Identical relations of exterior differential forms are a mathematical
expression of various types of conjugacy that leads to closed exterior
forms.

Such relations like
the Poincare invariant, vector and tensor identical relations,
the Cauchi-Riemann conditions, canonical relations, the integral
relations by Stokes, Gauss-Ostrogradskii, the thermodynamic
relations, the eikonal relations and so on are examples of identical relations
of closed exterior forms that have the form of relation (12) or its differential 
or integral analogs.

One can see that identical relations of closed exterior differential forms
make itself evident in various branches of physics and mathematics.

Below the mathematical and physical meaning of these relations and their role
in generating differential-geometrical structures will 
be disclosed with the help of evolutionary forms.

Thus, one can see that the closure conditions of exterior inexact
form and of corresponding dual form are a mathematical expression
of the conservation law and differential-geometrical structures.

And here there arise the questions of: (a) how closed inexact exterior
forms, which correspond to differential-geometrical structures and
reflect the properties of conservation laws, are obtained; (b) what
generates differential-geometrical structures;
and (c) what is responsible for such processes?

The mathematical apparatus of evolutionary differential forms enables
us to answer these questions.

\bigskip
{\bf Properties  of evolutionary forms.}

Above it has been shown that the evolutionary form commutator
includes the commutator of the manifold metric form which is nonzero.
Therefore, the evolutionary form commutator cannot be equal to zero.
This means that the evolutionary form
differential is nonzero. Hence, the evolutionary form, in
contrast to the case of the exterior form, cannot be closed. This leads to
that in the mathematical apparatus of evolutionary forms there arise
new unconventional elements like nonidentical relations and degenerate
transformations. Just such peculiarities allow to describe
evolutionary processes.

Nonidentical relations  can be written as
$$
d\phi=\eta^p\eqno(13)
$$
Here $\eta^p$ is the $p$-degree evolutionary form that is
unclosed, $\phi$ is some form of degree $(p-1)$, and
the differential $d\phi$ is a closed form of degree $p$.

In the left-hand side of this relation it stands the form differential,
i.e. a closed form that is an invariant object. In the right-hand
side it stands the unclosed form that is not an invariant
object. Such a relation cannot be identical one.

One can see a difference of relations for exterior forms and evolutionary
ones. In the right-hand side of identical relation (see relation (12))
it stands a closed form, whereas the form in the right-hand side of
nonidentical relation (13) is an unclosed one.

Nonidentical relations are obtained while describing any processes.
A relation of such type is obtained while analyzing the integrability
of the partial differential equation. An equation is integrable
if it can be reduced to the form $d\phi=dU$. However it
appears that, if the equation is not subject to an additional
condition (the integrability condition), it is reduced to the
form (13), where $\eta^p$ is an unclosed form and it cannot be
expressed as a differential.

Nonidentical relations of evolutionary forms are evolutionary relations
because they include the evolutionary form.
Such nonidentical evolutionary relations appear to be selfvarying
ones.  A variation of any object of the relation in some
process leads to a variation of another object and, in turn, a variation
of the latter leads to a variation of the former. Since one of the objects
is a noninvariant (i.e. unmeasurable)
quantity, the other cannot be compared with the first one, and hence,
the process of mutual variation cannot be completed.

The nonidentity of evolutionary relation is connected with a
nonclosure of the evolutionary form, that is, it is connected with
the fact that the evolutionary form commutator is nonzero.
As it has been pointed out, the evolutionary form commutator includes two terms.
one term specifies the mutual variations of the evolutionary form
coefficients, and the second term (the metric form commutator) specifies
the manifold deformation. These terms have a different nature and cannot
make the commutator vanish. In the process of selfvariation of the nonidentical
evolutionary relation it proceeds an exchange between the terms of the evolutionary
relation and this is realized according to the evolutionary relation. The evolutionary
form commutator describes a quantity that is a moving force of the evolutionary
process and leads to generation of differential-geometrical structures.

The significance of the evolutionary relation selfvariation consists in
the fact that in such a process it can be realized conditions under
which the identical relation is obtained from the nonidentical relation.
These are conditions of a degenerate transformation. The point of time at 
which such conditions are realized is that of originating the element of the
differential-geometrical structure.

\bigskip
{\bf Origination of differential-geometrical structures.}

To obtain the differential-geometrical structure, it is necessary
to obtain a closed inexact exterior form, i.e. the form closed on
pseudostructure. The condition of realization of the
pseudostructure is vanishing the interior differential of the
metric form. This is the closure condition for dual form and this leads
to realization of closed inexact exterior form.

Since the evolutionary form differential is nonzero, whereas the closed
exterior form differential is zero, the transition from the evolutionary
form to the closed exterior form is allowed only as a degenerate
transformation, i.e. a transformation that does not conserve the differential
The conditions of vanishing interior differential of the metric form
(the additional condition) are the conditions of degenerate
transformation.

As it has been already mentioned, the evolutionary differential form
$\eta^p$ involved into nonidentical relation (13) is an unclosed one.
The commutator, and hence the differential, of this form is nonzero.
That is,
$$d\eta^p\ne 0 \eqno(14)$$
If the conditions of degenerate transformation are realized, then from
the unclosed evolutionary form one can obtain a differential form closed
on pseudostructure. The differential of this form equals zero. That is,
it is realized the transition

$d\eta^p\ne 0 \to $ (degenerate transformation)
$\to d_\pi{}^*\eta^p=0$, $d_\pi \eta^p=0$.

The relations obtained
$$d_\pi \eta^p=0,  d_\pi{}^*\eta^p=0 \eqno(15)$$
are the closure conditions for exterior inexact form, and this points to 
realization of exterior form closed on pseudostructure, that is, this
points to origination of the differential-geometrical structure.

The conditions of degenerate transformation that lead to
origination of the differential-geometrical structures can be
connected with any symmetries.  (While describing material system, the 
symmetries can be conditioned, for example, by degrees of freedom
of material system). Since the conditions of degenerate
transformation are those of vanishing the interior differential of
metric form, that is, vanishing the interior (rather then total)
metric form commutator, the conditions of degenerate
transformation can be caused by symmetries of coefficients of the
metric form commutator (for example, it can be the symmetric 
connectedness).

To the conditions of degenerate transformation there
correspond a requirement that some functional expressions become equal
to zero. Such functional expressions are Jacobians, determinants,
the Poisson brackets, residues, and others.

Under the degenerate transformation on pseudostructure the evolutionary
form commutator vanishes, and this corresponds to a realization of the closed
inexact exterior form. But in this case the total evolutionary form commutator
is nonzero.

Vanishing on pseudostructure the external form differential (that
is, vanishing on pseudostructure the interior commutator of the
evolutionary form) points to that the exterior unclosed form is a
conservative quantity in the direction of pseudostructure.
However, in the direction normal to pseudostructure this quantity
exhibits a discontinuity. The value of such discontinuity is defined by 
the value of the evolutionary form commutator being nonzero. This 
argues to discreteness of the differential-geometrical structures.

Thus, while selfvariation of the evolutionary nonidentical
relation the interior metric form commutator can vanish (due to the
symmetries of coefficients of the metric form commutator).
This means that it is formed the pseudostructure on which
the differential form turns out to be closed. The emergency
of the form being closed on pseudostructure points out to origination
of the differential-geometrical structures.

It has been already noted that the conditions of degenerate transformations,
which lead to emergency of differential-geometrical structures,
can just be realized under selfvariation of the nonidentical
evolutionary relation.

\bigskip
{\it Obtaining identical relation from nonidentical one.}

On the pseudostructure $\pi$ evolutionary relation (13) converts to
the relation
$$
d_\pi\phi=\eta_\pi^p\eqno(16)
$$
which proves to be an identical relation. Indeed, since the form
$\eta_\pi^p$ is a closed one, on the pseudostructure this form turns
out to be a differential of some differential form. 
There are differentials in the left-hand and right-hand sides of
this relation. This means that the relation is an identical one.

From evolutionary relation (13) it is obtained the identical on the
pseudostructure relation. In this case the evolutionary relation itself
remains to be nonidentical one. (At this point it should be
emphasized that differential, which equals zero, is an interior one.
Under degenerate transformation the evolutionary form commutator  becomes
zero only on the pseudostructure. The total evolutionary form commutator
is nonzero. The evolutionary form remains to be unclosed.)

It can be shown that all identical relations of the exterior
differential form theory are obtained from nonidentical relations (that
contain the evolutionary forms) by applying degenerate transformations.

The degenerate transform is realized as a transition to
nonequivalent coordinate system: a transition from the 
noninertial coordinate system to the locally inertial that.
Evolutionary relation (13) and condition (14) relate to the system
being tied to the deforming manifold, whereas
condition (15) and identical relations (16) may relate only to the
locally inertial coordinate system being tied to a pseudostructure.

Transition from nonidentical relation (13) to identical relation (16)
means the following. Firstly, it is from such a relation that
one can obtain the differential $d_\pi\phi$ and find the desired
function $\phi_\pi$. And, secondly, an emergency
of the closed (on pseudostructure) inexact exterior form $\eta_\pi^p$
(right-hand side of relation (16)) points to an origination of the
conservative object - the differential-geometrical structure.
It occurs that the emergency of the differential-geometrical structure
is connected with the realization of the differential $d_\pi\phi$.
Below it will be shown that the function $\phi_\pi$ is the state-function
of the system described, which generates the differential-geometrical
structures.

Thus, the mathematical apparatus of evolutionary forms
describes the process of generation of the closed inexact exterior
forms, and this discloses the process of origination of
the differential-geometrical structure, namely, a new conjugated object.

The above described process of generating the differential-geometrical
structures discloses the process of conjugating any elements or operators.

The evolutionary differential form is an unclosed form, that is, it is
a form whose differential is not equal to zero. The differential of
the exterior differential form equals zero. To the closed exterior form
there correspond conjugated operators, whereas to the evolutionary form
there correspond nonconjugated operators. A transition from the
evolutionary form to the closed exterior form and origination of the
differential-geometrical structures is a transition from nonconjugated
operators to conjugated ones. This is expressed mathematically
as a transition from a nonzero differential (the evolutionary form
differential is nonzero) to a differential that equals zero (the closed
exterior form differential equals zero).

It can be seen that the process of generating the differential-geometrical
structures is a mutual exchange between the quantities of different
nature (for example, between the algebraic and geometric quantities,
between the physical and spatial quantities) and an exchange while
realizing any symmetry.

\bigskip
{\it Characteristics of differential-geometrical structures.}

Since the closed exterior form, which corresponds to the
differential-geometrical structure emerged, was obtained from
the evolutionary form, it is evident that characteristics
of this structure has to be connected with those of the evolutionary form
and of the manifold on which this form is defined, with the conditions
of degenerate transformation and with the values of commutators of the
evolutionary form and the manifold metric form.

The conditions of degenerate transformation, as it was said before,
determine the pseudostructures. The first term of the evolutionary form
commutator determines the value of the discrete change (the quantum),
which the quantity conserved on the pseudostructure undergoes when
transition from one pseudostructure to another. The second term of the
evolutionary form commutator specifies a characteristics that fixes the
character of the initial manifold deformation, which took place before
the differential-geometrical structure emerged.  (Spin is an example of such
a characteristics).

The connection of the differential-geometrical structures with the
skew-symmetric differential forms allows to introduce a classification 
of these structures in dependence on parameters that specify the
skew-symmetric differential forms and enter into nonidentical and
identical relation. To determine these parameters one has to consider
the problem of integration of the nonidentical evolutionary relation.

Under degenerate transformation from the nonidentical evolutionary
relation one obtains a relation being identical on pseudostructure.
Since the right-hand side of such a relation can be expressed in terms
of differential (as well as the left-hand side), one obtains a relation
that can be integrated, and as a result he obtains a relation with the
differential forms of less by one degree.

The relation obtained after integration proves to be nonidentical
as well.

By sequential integrating the nonidentical relation of degree $p$ (in
the case of realization of the corresponding degenerate transformations
and forming the identical relation), one can get closed (on the
pseudostructure) exterior forms of degree $k$, where $k$ ranges
from $p$ to $0$.

In this case one can see that after such integration the closed (on the
pseudostructure) exterior forms, which depend on two parameters, are
obtained. These parameters are the degree of evolutionary form $p$
in the evolutionary relation and the degree of created closed forms $k$.

In addition to these parameters, another parameter appears, namely, the
dimension of space. If the evolutionary relation generates the closed
forms of degrees $k=p$, $k=p-1$, \dots, $k=0$, to them there correspond
the pseudostructures of dimensions $(N-k)$, where $N$ is the space (formed) 
dimension. 

\bigskip
{\it Forming fields and manifolds.}

The pseudostructures, on which the closed {\it inexact} forms are
defined, form the pseudomanifolds. (Integral surfaces, pseudo-Riemann
and pseudo-Euclidean spaces are the examples of such manifolds). In this
process the dimensions of the  manifolds formed  are connected with the
evolutionary form degree.

To transition from pseudomanifolds to metric manifolds it is assigned
a transition from closed {\it inexact} differential forms to {\it exact}
exterior differential forms. (Euclidean and Riemann spaces are examples
of metric manifolds).

Since the closed metric form is dual with respect to some closed exterior
differential form, the metric forms cannot become closed by themselves,
independently of the exterior differential form. This proves that manifolds
with closed metric forms are connected with the closed exterior
differential forms. This indicates that the fields of conservative
quantities are formed from closed exterior forms at the same instant of 
time when the manifolds are created from the pseudoctructures. (The specific
feature of the manifolds with closed metric forms that have been formed
is that they can carry some information.) That is, the closed exterior
differential forms and manifolds, on which they are defined, are
mutually connected objects.

\bigskip
{\it Symmetries, conservation laws, differential-geometrical structures.}

As it has been noted, the closure conditions of exterior inexact
form and of corresponding dual form are a mathematical expression
of the conservation law and differential-geometrical structures.

The closure of the exterior differential forms and corresponding dual form
is the result of the conjugacy of elements of exterior or dual
forms.

Conjugacy is possible if there is one or other type of symmetry.

It has been shown that the differential-geometrical structures
originate when the conditions of degenerate transformation are realized.
These conditions relate to the symmetries of coefficients of metric form of
the manifold on which the evolutionary form is defined. These symmetries of
the dual form are exterior symmetries. Symmetries of the closed exterior form
that is formed on the pseudostructure are interior symmetry.

Whereas the exterior symmetries of the dual forms
are connected with degenerate transformations, the interior symmetries, i.e.
those of the closed exterior forms, are connected with nonidentical
transformations.

Degenerate and nondegenerate transformations are mutually
connected. Degenerate transformations lead to formatting the
differential-geometrical structures and nondegenerate
transformations execute a transition from one
differential-geometrical structure to another.

One of the problems in the theory of symmetry is a searching for symmetries of
differential equations. A knowledge of symmetries enables one to get
a solution of differential equations that corresponds to the conservation law
and defines the differential-geometrical structures. The dependence of
symmetries on any parameter, i.e. the dependence of the nondegenerate transformations
on a given parameter allows to study the ``evolution" of the
differential-geometrical structures in given parameter. But in this case the question
of how these structures emerge is not posed (that is, the evolutionary process 
of originating these structures is not considered). As it has been shown in present 
paper, the answer to this question gives the theory of evolutionary forms.

It has been already noted that the evolutionary forms which 
generate the closed exterior forms corresponding to the 
differential-geometrical structures appear while describing any 
processes by differential equations. It can be shown that the tangent 
manifold of differential equation is nondifferentiable manifold (the
Lagrangian manifold is an example) and the derivatives satisfying
the differential equation do not make up a differential 
(closed exterior forms). However, under degenerate
transformation it occurs a transition from the tangent manifold
to cotangent one (the manifold of Hamiltonian systems is an
example) on which the differential made up by derivatives of
differential equation becomes a closed form. This points to originating 
differential-geometrical structures. That is, the origination of the 
differential-geometrical structure is connected with  degenerate 
transformation that 
executes the transition from tangent space to cotangent one. And
nondegenerate transformation executes the transition into tangent
space from any differential-geometrical structure to another. 

A peculiarity of the degenerate transformation can be considered by the
example of Hamiltonian systems.
Here the degenerate transformation is a transition from the tangent
space ($q_j,\,\dot q_j)$), which is a tangent (Lagrangian) manifold,
to the cotangent characteristical (Hamiltonian) manifold
($q_j,\,p_j$). On the other hand, the nondegenerate canonical
transformation is a transition from one characteristical manifold
($q_j,\,p_j$) to another characteristical manifold ($Q_j,\,P_j$).
$\{$The formula of nondegenerate canonical
transformation can be written as $p_jdq_j=P_jdQ_j+dW$,
where $W$ is the generating function$\}$.

\bigskip
{\bf Physical meaning of the differential-geometrical structures.}

The differential-geometrical structures obtained may carry a physical
meaning.

It has been already noted that the evolutionary forms appear under description
of a certain process. In particular, they are obtained while describing physical
processes in material media (material systems). Analysis of the evolutionary forms
obtained and nonidentical relations shows that material media generate
differential-geometrical structures, which are physical structures forming
physical fields.
\{Material system is a variety of elements that have internal structure
and interact to one another. As examples of material systems it may be
thermodynamic, gas dynamical, cosmic systems, systems of elementary
particles and others.\}

Evolutionary forms, as well as the closed exterior forms, correspond to conservation
laws. But closed exterior forms correspond to conservation laws that can be named as
exact ones, whereas evolutionary forms correspond to the balance conservation laws.

\bigskip
{\it The balance conservation laws.}

The balance conservation laws are conservation laws that
establish a balance between the variation of physical quantity and
the corresponding external action. They are described by differential equations.
The balance conservation laws for material systems are conservation laws
for energy, linear momentum, angular momentum, and
mass. From the equations, which describe the balance conservation laws,
one can obtain a relation that is nonidentical relation since it contains 
the evolutionary forms.

Let us analyze the equations that describe the balance conservation
laws for energy and linear momentum.

In the accompanying reference system (this system is connected with the manifold
built by the trajectories of the material system elements) the energy
equation is written in the form
$$
{{\partial \psi }\over {\partial \xi ^1}}\,=\,A_1 \eqno(17)
$$

Here $\psi$  is the functional specifying the state of material system
(the action functional, entropy, wave function can be regarded as
examples of the functional), $\xi^1$ is the coordinate along the
trajectory, $A_1$ is the quantity that depends on specific features of
material system and on external energy actions onto the system.

In a similar manner, in the accompanying reference system the
equation for linear momentum appears to be reduced to the equation of
the form
$$
{{\partial \psi}\over {\partial \xi^{\nu }}}\,=\,A_{\nu },\quad \nu \,=\,2,\,...\eqno(18)
$$
where $\xi ^{\nu }$ are the coordinates in the direction normal to the
trajectory, $A_{\nu }$ are the quantities that depend on the specific
features of material system and on external force actions.

Eqs. (17) and (18) can be convoluted into the relation
$$
d\psi\,=\,A_{\mu }\,d\xi ^{\mu },\quad (\mu\,=\,1,\,\nu )\eqno(19)
$$
where $d\psi $ is the differential
expression $d\psi\,=\,(\partial \psi /\partial \xi ^{\mu })d\xi ^{\mu }$.

Relation (19) can be written as
$$
d\psi \,=\,\omega \eqno(20)
$$
here $\omega \,=\,A_{\mu }\,d\xi ^{\mu }$ is the skew-symmetrical differential form of the first degree.

Since the balance conservation laws are evolutionary ones, the relation
obtained is also an evolutionary relation.

Relation (20) was obtained from the equation of the balance conservation
laws for energy and linear momentum. In this relation the form $\omega $
is that of the first degree. If the equations of the balance conservation
laws for angular momentum be added to the equations for energy and linear
momentum, this form in the evolutionary relation will be the form of the
second degree. And in  combination with the equation of the balance
conservation law of mass this form will be the form of degree 3.

Thus, in general case the evolutionary relation can be written as
$$
d\psi \,=\,\omega^p \eqno(21)
$$
where the form degree  $p$ takes the values $p\,=\,0,1,2,3$..
(The evolutionary relation for $p\,=\,0$ is similar to that in the differential forms, and it was
obtained from the interaction of energy and time.)

Let us show that relation obtained from the equation
of the balance conservation laws proves to be nonidentical.

To do so we shall analyze relation (20).

In the left-hand side of relation (20) there is a
differential that is a closed form. This form is an invariant
object. The right-hand side of relation (20) involves the differential
form $\omega$, that is not an invariant object because in real processes,
as it is shown below, this form proves to be unclosed. The commutator of this
form is nonzero. The components of commutator of the form $\omega \,=\,A_{\mu }d\xi ^{\mu }$
can be written as follows:
$$
K_{\alpha \beta }\,=\,\left ({{\partial A_{\beta }}\over {\partial \xi ^{\alpha }}}\,-\,
{{\partial A_{\alpha }}\over {\partial \xi ^{\beta }}}\right )
$$
(here the term  connected with the manifold metric form 
has not yet been taken into account).

The coefficients $A_{\mu }$ of the form $\omega $ have been obtained either
from the equation of the balance conservation law for energy or from that for
linear momentum. This means that in the first case the coefficients depend
on the energetic action and in the second case they depend on the force action.
In actual processes energetic and force actions have different nature and appear
to be inconsistent. The commutator of the form $\omega $ consisted of 
the derivatives of such coefficients is nonzero.
This means that the differential of the form $\omega $
is nonzero as well. Thus, the form $\omega$ proves to be unclosed and
cannot be a differential like the left-hand side.

This means that the relation (20) cannot be an identical one.

What is a physical meaning of such a relation?

This relation obtained from the equations of the balance
conservation laws involves the functional that specifies the
material system state. However, since this relation turns out to
be not identical, from this relation one cannot get the
differential $d\psi $  that could point out to the equilibrium state of
material system. The absence of differential means that the system
state is nonequilibrium. That is, in material system the internal
force acts. This leads to distortion of trajectories of material
system. A manifold made up by the trajectories (the accompanying
manifold) turns out to be a deforming manifold. The differential
form $\omega$,  as well as the forms $\omega^p$, appear to be
evolutionary forms. Commutators of these forms will contain an additional term
connected with the commutator of unclosed metric form of manifold.

The availability of two terms in the commutator of the form $\omega^p $,
as it has been already shown, leads to that the nonidentical evolutionary
relation turns out to be a selfvarying relation.

Selfvariation of the nonidentical evolutionary relation points to the fact
that the nonequilibrium state of material system turns out to be selfvarying.
It is evident that this selfvariation proceeds under the action of internal force
whose quantity is described by commutator of the unclosed evolutionary form
$\omega^p $. (If the commutator
be zero, the evolutionary relation would be identical, and this would
point to the equilibrium state, i.e. the absence of internal forces.)
Everything that gives a contribution into the commutator of the form
$\omega^p $ leads to emergency of internal force.

Above it has been shown that under degenerate transformation from nonidentical
evolutionary relation it can be obtained the identical relation
$$
d_\pi\psi=\omega_\pi^p\eqno(22)
$$

From such a relation one can obtain the state function and this
corresponds to equilibrium state of the system. But identical
relation can be realized only on pseudostructure (which is
specified by the condition of degenerate transformation). This
means that a transition of material system to equilibrium state
proceeds only locally. In other words, it is realized a transition
of material system from nonequilibrium state to locally
equilibrium one. In this case the general state of material system
remains to be nonequilibrium.

It has been already noted that the symmetries of coefficients of the metric form
commutator are conditions of degenerate transformation. They are conditioned by
degrees of freedom of material system.

As one can see from the analysis of nonidentical evolutionary
relation, the transition of material system from nonequilibrium
state to locally-equilibrium state proceeds spontaneously in the
process of selfvarying nonequilibrium state of material system
under realization of any degrees of freedom of this system.
(Translational degrees of freedom, internal degrees of freedom of
the system elements, and so on can be examples of such degrees of
freedom).

As it has been already said above, the transition from nonidentical relation (21)
obtained from the balance conservation laws to identical relation (22) means
the following.
Firstly, an existence of the state differential (left-hand side
of relation (22)) points to a transition of the material system from nonequilibrium state
to the locally-equilibrium state. And, secondly, an emergency of the
closed (on pseudostructure) inexact exterior form (right-hand side
of relation (22)) points to an origination of the physical structure.

Thus one can see that the transition of material system from
nonequilibrium state to locally-equilibrium state is accompanied
by originating differential-geometrical structures, which are
physical structures. The emergency of physical structures in the
evolutionary process reveals in material system as an emergency of
certain observable formations, which develop spontaneously. In
this manner the causality of emerging various observable
formations in material media is explained. Such formations and
their manifestations are fluctuations, turbulent pulsations,
waves, vortices, creating massless particles and others. The
intensity of such formations is controlled by a quantity
accumulated by the evolutionary form commutator at the instant in
time of originating physical structures.

Physical structures that are generated by material systems made up
physical fields.

The availability of physical structures points out to the fulfilment
of conservation laws. These are conservation laws for physical fields.
The process of generating physical fields demonstrates a connection of
these conservation laws, which can be referred to as exact ones, with the balance
(differential) conservation laws for material media. 
Closed inexact exterior
forms that correspond to physical structures and exact conservation laws
are obtained from equations describing the balance conservation laws.

Since the closed exterior forms corresponding to physical structures are
obtained from the evolutionary forms describing material systems, the
characteristics of physical structures are determined by characteristics of
the material system generating these structures, and the parameters of
evolutionary forms and closed exterior forms enables one to classify
physical structures.

As it has been shown above, the type of differential-geometrical structures, and 
hence of physical structures (and,
accordingly, of physical fields) generated by the evolutionary
relation, depends on the degree of differential forms $p$ and $k$
and on the dimension of original inertial space $n$ (here $p$ is
the degree of the evolutionary form in the nonidentical relation
that is connected with a number of interacting balance
conservation laws, and $k$ is the  degree of closed form
generated by the nonidentical relation). By introducing the
classification by numbers $p$, $k$, $n$ one can understand the
internal connection between various physical fields. Since physical
fields are the carriers of interactions, such
classification enables one to see a connection between
interactions.  It can be shown that the case $k=0$ corresponds to 
strong interaction, $k=1$ corresponds to weak interaction,
$k=2$ corresponds to electromagnetic interaction, and $k=3$ corresponds
to gravitational interaction.

Mathematical apparatus of closed exterior forms that corresponds to physical
structures and conservation laws for physical fields lies at the basis of
field theory.

It can be shown that the equations of existing field theories are
those obtained on the basis of the properties of the exterior form
theory. The Hamilton formalism is based on the properties of closed
exterior and dual forms of the first degree, quantum mechanics does on the forms
of zero degree, the electromagnetic field equations are based on the forms of
second degree. The third degree forms are assigned to the gravitational field.

The gauge symmetries, which are interior symmetries of the field
theory equations, are symmetries of closed (inexact) forms.

The internal symmetries in field theory are those of closed exterior
differential forms, whereas the external symmetries in field theory
are symmetries of relevant dual forms.

The gauge transformations of field theory, which are nondegenerate
transformations of the closed exterior differential forms, are connected with
the gauge symmetries. Since the closed exterior differential form is a differential
(a total one if the form is exact, or an interior one on pseudostructure
if the form is inexact), it is
obvious that the closed form proves to be invariant under all
transformations that conserve the differential. The unitary transformations
(0-form), the tangent and canonical transformations (1-form), the gradient and
gauge transformations (2-form) and so on are  examples of such transformations.
{\it These are gauge transformations for spinor, scalar, vector, and tensor
(3-form) fields}.

It has been shown that the closed exterior forms and relevant dual forms, which
correspond to the conservation laws for physical fields, are obtained from
the evolutionary forms, which describe the balance conservation laws
for material media. This proceeds under degenerate transformation, which is
connected with the degrees of freedom of material system. The conditions of
degenerate transformation defines a closed dual form. From this it follows
that the external symmetries, namely, the symmetries of dual forms, are due to
the degrees of freedom of material system. It is for this reason the exterior
symmetries are spatial-temporal symmetries.

The realization of the closed dual form,
which proceeds due to the degrees of freedom of material system,
leads to realization of the closed exterior form, that is, to the
conjugacy  of the differential form elements, and emergency of
internal symmetries. From this one can see a connection between
internal and external symmetries.

Whereas the internal symmetries are connected with the conservation laws for
physical fields, the external symmetries caused by the degrees of freedom of
material media are connected with the balance conservation laws for material
media.

The nondegenerate transformations are connected with internal symmetries,
and the degenerate transformations of evolutionary forms
are connected with external symmetries.

\bigskip
\centerline {\bf Conclusion}

As it has been shown the theory of evolutionary forms explains
the process of generation of differential-geometrical structures.
This cannot be carried out within the framework of any other formalisms.

1. Schutz B.~F., Geometrical Methods of Mathematical Physics. Cambrige 
University Press, Cambrige, 1982.

2. Bott R., Tu L.~W., Differential Forms in Algebraic Topology. 
Springer, NY, 1982.

3. Encyclopedia of Mathematics. -Moscow, Sov.~Encyc., 1979 (in Russian).

\end{document}